\documentclass{amsart}
\usepackage{latexsym,amssymb,amsfonts,amscd,emlines2,amsthm,amstext}
\unitlength=1mm\linethickness{0.6pt}
\numberwithin{equation}{section}
\newtheorem{conj}{Conjecture}[section]
\newtheorem{antipode}{Proposition}[section]
\newtheorem{braided}{Proposition}[section]

           \theoremstyle{definition}
\newtheorem{braided coalgebra}[Yetter]{Definition}

\begin{document}
\title[Clifford Hopf-gebra]{Clifford Hopf-gebra and Bi-universal
Hopf-gebra}\author[Zbigniew Oziewicz]{Zbigniew Oziewicz}
\address{University of Wroc{\l}aw, Institute of Theoretical Physics,
plac Maksa Borna 9, 50204  Wroc{\l}aw, Poland}
\curraddr{Universidad Nacional Aut\'onoma de M\'exico, Facultad de
Estudios Superiores,  C.P. 54700 Cuautitl\'an Izcalli, Apartado Postal
\# 25, Estado de M\'exico}
\email{oziewicz@servidor.unam.mx, oziewicz@ift.uni.wroc.pl}
\subjclass{Primary 16W30; 17A42}
\keywords{Clifford Hopf-gebra, Clifford bi-gebra, $n$-ary bigebra,
braided bi-gebra}\thanks{The Author is a member of Sistema Nacional de
Investigadores, M\'exico.}
\thanks{Submitted August 28, 1997, to Czechoslovak Journal of Physics,
q-alg/9709016. This paper is in final form and no version of it will be
submitted for publication elsewhere.}

\newcommand{\N}{\mathbb{N}}\newcommand{\C}{\mathbb{C}}
\newcommand{\R}{\mathbb{R}}\newcommand{\Z}{\mathbb{Z}}
\newcommand{\B}{\mathcal{B}}
\newcommand{\ie}{\textit{i.e.}$\,$}
\newcommand{\be}{\begin{equation}}\newcommand{\ee}{\end{equation}}
\newcommand{\ba}{\begin{array}}\newcommand{\ea}{\end{array}}
\newcommand{\ra}{\longrightarrow}
\newcommand{\id}{\operatorname{id}}
\newcommand{\der}{\text{der}}
\newcommand{\alg}{\text{alg}}\newcommand{\cog}{\text{cog}}
\newcommand{\im}{\text{im}}
\newcommand{\End}{\operatorname{End}}
\newcommand{\Nat}{\operatorname{Nat}}
\newcommand{\gen}{\text{gen}}
\newcommand{\lin}{\operatorname{lin}}\newcommand{\sh}{\text{sh}}
\newcommand{\hopf}{\text{hopf}}
\newcommand{\Span}{\text{span}}
\newcommand{\grade}{\text{grade}}
\newcommand{\spec}{\text{spec}}
\newcommand{\eid}{\operatorname{id}}
\newcommand{\idL}{{\id}_{L^{\otimes 2}}}
\newcommand{\eidL}{\operatorname{id_{L^{\otimes 2}}}}
\newcommand{\eder}{\operatorname{der}}
\newcommand{\ering}{\operatorname{ring}}
\newcommand{\emod}{\operatorname{mod}}
\newcommand{\ealg}{\operatorname{alg}}
\newcommand{\eim}{\operatorname{im}}
\newcommand{\eEnd}{\operatorname{End}}
\newcommand{\egen}{\operatorname{gen}}
\newcommand{\elin}{\operatorname{lin}}
\newcommand{\edeg}{\operatorname{deg}}
\newcommand{\edim}{\operatorname{dim}}
\newcommand{\egrade}{\operatorname{grade}}
\newcommand{\Cl}{\mathcal{C}\ell}
\newcommand{\half}{\textstyle{\frac{1}{2}}}

\begin{abstract} For a scalar product $\xi$ on co-vectors, the Clifford
co-product $\triangle^\xi$ of multivectors is calculated from the dual
Clifford algebra. With respect to this co-product $\triangle^\xi,$ unit
is not group-like and vectors are not primitive. For a scalar product
$\eta$ on vectors the Clifford product $\wedge^\eta$ and the Clifford
co-product $\triangle^\xi$ fits to the bi-gebra with respect to the
family of the (pre)-braids. The Clifford bi-gebra is in a braided
category iff $\xi=0$ or $\eta=0.$\end{abstract}\maketitle
\tableofcontents

\section{Multi-ary Bi-gebra} A multiplicative category is a category
with a bifunctor of bin-ary operation, an anihilation $2\rightarrow 1$,
denoted by two initial leaves and one node. A co-multiplicative
category possess a binary co-operation, a creation process
$1\rightarrow 2$. A pairing $2\rightarrow 0,$ a bin-ary anihilation
$2\rightarrow 1,$ a bin-ary creation $1\rightarrow 2,$ and bin-ary
scattering $2\rightarrow 2$ are represented by the prime graph nodes in
Diagram 1. All diagrams are directed and is recommendent to read them
from the top to the bottom.

$$\begin{picture}(65,15)          
\put(8,15){\oval(6,20)[b]}
\put(23,15){\oval(6,10)[b]}
\put(38,5){\oval(6,10)[t]}
\emline{38}{10}{1}{38}{15}{2}
\emline{23}{10}{3}{23}{5}{4}
\emline{50}{15}{5}{60}{5}{6}
\emline{60}{15}{7}{56}{11}{8}
\emline{54}{9}{9}{50}{5}{10}
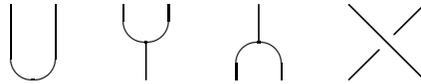
\put(32,3){\makebox(0,0)[ct]{\begin{tabular}{c}\textsc{Diagram} 1.
A pairing, binary multiplication (anihilation),\\\hspace{2cm} binary
co-product \ie creation and a scattering (braid)\end{tabular}}}
\end{picture}$$\smallskip

Throughout this paper $\Bbbk$ denotes a commutative ring.
A pair of $\Bbbk$-modules, say $A$ and $B,$ with a pairing $A\otimes
B\rightarrow\Bbbk,$ is said to be a dual pair. Dual pair of al-gebras
(co-gebras) extends to pair of co-gebras (al-gebras) and these
structures, al-gebra and co-gebra, may close to bi-gebra for a family
of pre-braids. We calculate these structures from the assumptions
displayed on Diagrams 2-4.

$$\begin{picture}(100,25)         
\put(12.50,15){\oval(15,20)[b]}
\put(12.50,15){\oval(5,10)[b]}
\put(7.50,15){\oval(5,10)[t]}
\emline{7.33}{20}{1}{7.33}{23}{2}
\emline{15}{23}{3}{15}{15}{4}
\emline{20}{15}{5}{20}{23}{6}
\put(25,14){\makebox(0,0)[cc]{$=$}}
\put(38,23){\oval(6,18)[b]}
\put(34,14){\oval(8,18)[b]}
\emline{30}{23}{7}{30}{14}{8}
\put(68.50,15){\oval(15,20)[b]}
\put(68.50,15){\oval(5,10)[b]}
\put(73.50,15){\oval(5,10)[t]}
\emline{73.33}{20}{9}{73.33}{23}{10}
\emline{66}{23}{11}{66}{15}{12}
\emline{61}{15}{13}{61}{23}{14}
\put(81,14){\makebox(0,0)[cc]{$=$}}
\put(89,23){\oval(6,18)[b]}
\put(93,14){\oval(8,18)[b]}
\emline{97}{23}{15}{97}{14}{16}
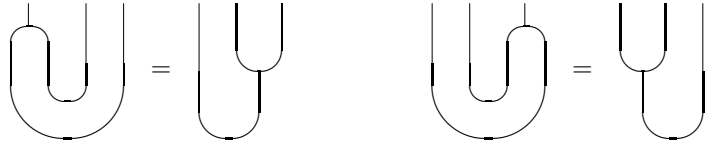
\put(50,0){\makebox(0,0)[c,t]{\textsc{Diagram} 2. The product -
co-product duality}}\end{picture}$$

The Diagram 2 imply that whenever product (co-product) is in variety
then co-product (product) is in `the same' co-variety.
Diagram 2 extends to duality between $n$-ary multiplication and
$n$-ary co-multiplication as shown for tern-ary operations,
$1\leftrightarrow 3,$ on Diagram 3.

$$\begin{picture}(54,30)       
\put(8,14){\oval(6,10)[t]}
\put(12.50,23){\oval(9,30)[b]}
\put(12.50,14){\oval(15,18)[b]}
\put(12.50,14){\oval(3,6)[b]}
\emline{14}{14}{1}{14}{23}{2}
\emline{20}{23}{3}{20}{14}{4}
\put(30,14){\makebox(0,0)[cc]{$=$}}
\put(43,23){\oval(6,36)[b]}
\put(46,23){\oval(6,18)[b]}
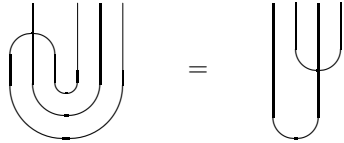
\put(27,3){\makebox(0,0)[ct]{\begin{tabular}{c}\textsc{Diagram} 3.
The tern-ary co-product - product duality\end{tabular}}}\end{picture}$$

$$\begin{picture}(43,30)                
\put(8,8.17){\oval(6,15.67)[t]}
\put(8,30){\oval(6,16)[b]}
\put(16,19){\makebox(0,0)[cc]{$=$}}
\put(24,22){\oval(6,10)[t]}
\put(35,22){\oval(6,10)[t]}
\put(24,16){\oval(6,10)[b]}
\put(35,16){\oval(6,10)[b]}
\emline{8}{22}{1}{8}{16}{2}
\emline{24}{30}{3}{24}{27}{4}
\emline{24}{11}{5}{24}{8}{6}
\emline{35}{8}{7}{35}{11}{8}
\emline{35}{27}{9}{35}{30}{10}
\emline{21}{22}{11}{21}{16}{12}
\emline{38}{16}{13}{38}{22}{14}
\emline{27}{22}{15}{32}{16}{16}
\emline{32}{22}{17}{30}{20}{18}
\emline{29}{18}{19}{27}{16}{20}
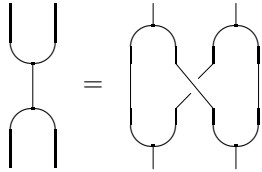
\put(20,4){\makebox(0,0)[c,t]{\textsc{Diagram} 4. The bin-ary bi-gebra
with one (pre)-braid}}\end{picture}$$

A tern-ary bi-gebra is defined on Diagram 5 and this made clear the
definition of the multi-ary bigebra.

$$\begin{picture}(60,35)    
\put(8,35){\oval(6,20)[b]}
\put(8,9){\oval(6,18)[t]}
\emline{8}{35}{1}{8}{9}{2}
\put(20,22){\makebox(0,0)[cc]{$=$}}
\put(33,18){\oval(6,10)[b]}
\put(43,18){\oval(6,10)[b]}
\put(53,18){\oval(6,10)[b]}
\put(33,25){\oval(6,10)[t]}
\put(43,25){\oval(6,10)[t]}
\put(53,25){\oval(6,10)[t]}
\emline{30}{25}{3}{30}{18}{4}
\emline{56}{18}{5}{56}{25}{6}
\emline{33}{35}{7}{33}{25}{8}
\emline{33}{25}{9}{40}{18}{10}
\emline{36}{25}{11}{50}{18}{12}
\emline{46}{25}{13}{53}{18}{14}
\emline{53}{18}{15}{53}{9}{16}
\emline{53}{35}{17}{53}{25}{18}
\emline{53}{25}{19}{50}{22}{20}
\emline{49}{21}{21}{48}{20}{22}
\emline{46.67}{18.67}{23}{46}{18}{24}
\emline{43}{35}{25}{43}{22.67}{26}
\emline{43}{20.33}{27}{43}{9}{28}
\emline{40}{25}{29}{39}{24}{30}
\emline{38.33}{23.33}{31}{37}{22}{32}
\emline{36}{21}{33}{33}{18}{34}
\emline{33}{18}{35}{33}{9}{36}
\emline{36}{18}{37}{38}{19}{38}
\emline{39.33}{19.67}{39}{42}{21}{40}
\emline{44}{22}{41}{46}{23}{42}
\emline{48}{24}{43}{50}{25}{44}
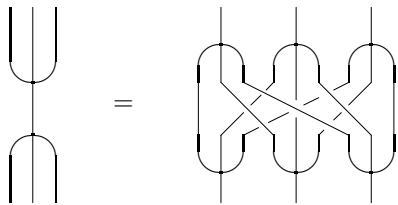
\put(30,7){\makebox(0,0)[ct]{\begin{tabular}{c}\textsc{Diagram} 5.
A tern-ary bi-gebra with nine scatterings (braids)\end{tabular}}}
\end{picture}$$

\section{Clifford Co-gebra and Antipode} Dirac in 1928 predicted the
existance of an anti-matter, spin $\half$ positrons, in terms of the
Clifford algebra. To understand an anti-matter for any spin we need an
action of the Clifford algebra $\Cl$ in a tensor product of
$\Cl$-modules. The differential Dirac operator for mesons of zero and
higher spins needs an action of the Clifford algebra on a tensor
product of Clifford algebras and this action is illustrated on Diagram
6. This action depends on choosen co-product, the simplest one is known
as the Duffin \& Kemmer \& Petiau co-product [Duffin 1938, Kemmer 1939,
1943]. The main problem is the classification of the co-products which
fit the Clifford algebra into bi-gebra. The present paper is the
introduction into this subject.

A category is said to be autonomous if
$\forall\;M\in\text{obj},$ $\exists\,!$ a left dual $M^*$ and a right
dual $^*M$ [Freyd \& Yetter 1992]. An autonomous category is said to
be pivotal if $M^*\simeq {^*M}.$

Let $M$ be $\Bbbk$-module and $\eta$ and $\xi$ be scalar products,
\begin{gather*}\eta\in\lin(M,M^*)\simeq M^{*\otimes 2}\simeq
\lin(M^{\otimes 2},\Bbbk),\\
\xi\in\lin(M^*,M)\simeq M^{\otimes 2}\simeq\lin(M^{*\otimes 2},\Bbbk).
\end{gather*}
A pair of the mutually dual Clifford $\Bbbk$-algebras is paired by
determinant (scalar product independent),
\begin{gather}\Cl(M,\eta)\simeq\{M^{\wedge},\wedge^\eta\},\quad
\Cl(M^*,\xi)\simeq\{M^{*\wedge},\wedge^\xi\},\label{clifford}\\
M^{*\wedge}\,\otimes\,M^{\wedge}\quad
\stackrel{\text{det}}{\longrightarrow}\quad\Bbbk.\end{gather}
A $\xi$-dependent co-multiplication
$\triangle^\xi:M^\wedge\longrightarrow M^\wedge\otimes M^\wedge$
is calculated from the product - co-product duality of Diagram 2.

$$\begin{picture}(30,35)      
\put(8,17){\oval(6,12)[b]}
\put(8,22){\oval(6,12)[t]}
\put(19,17){\oval(6,12)[b]}
\emline{11}{22}{1}{16}{17}{2}
\emline{22}{17}{3}{22}{31}{4}
\emline{16}{31}{5}{16}{22}{6}
\emline{16}{22}{7}{14.33}{20.33}{8}
\emline{12.67}{18.67}{9}{11}{17}{10}
\emline{11}{17}{11}{11}{16}{12}
\emline{5}{16}{13}{5}{22}{14}
\emline{8}{28}{15}{8}{31}{16}
\emline{8}{11}{17}{16}{7}{18}
\emline{16}{7}{19}{16}{5}{20}
\emline{22}{5}{21}{22}{8}{22}
\emline{22}{8}{23}{19}{11}{24}
\put(19,11){\circle*{1}}
\put(8,11){\circle*{1}}
\put(8,28){\circle*{1}}
\put(6,28){\makebox(0,0)[rb]{$\triangle^\xi$}}
\put(16,20){\makebox(0,0)[lc]{$\sigma$}}
\put(22,11){\makebox(0,0)[lc]{$\wedge^\eta\otimes\wedge^\eta$}}
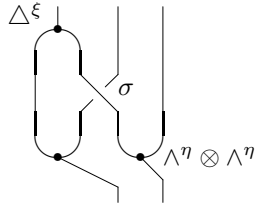
\put(15,2){\makebox(0,0)[ct]{\begin{tabular}{c}\textsc{Diagram} 6.
A co-product dependent action of $\Cl$\\ on tensor product of
$\Cl$-modules\end{tabular}}}\end{picture}$$\smallskip

In the sequel $\{e_\mu\in M\}$ denotes a basis and
$\{\varepsilon^\mu\in M^*\equiv\lin(M,\Bbbk)\}$ is a dual basis,
$\varepsilon^\mu e_\nu=\delta^\mu_\nu.$ For $1\in\Bbbk<M^\wedge$ and
$v,w\in M,$

\begin{multline*}\triangle^\xi 1=1\otimes 1+\sum\xi(\varepsilon^\mu
\otimes\varepsilon^\nu)e_\mu\otimes e_\nu\\
-\sum\xi(\varepsilon^{\mu_1}\wedge\varepsilon^{\mu_2},
\varepsilon^{\nu_1}\wedge\varepsilon^{\nu_2})
(e_{\mu_1}\wedge e_{\mu_2})\otimes(e_{\nu_1}\wedge e_{\nu_2})\\
-\sum\xi(\text{tri-co-vectors}^{\otimes 2})\text{trivectors}^{\otimes
2}+\ldots+(-)^{\left[\half\text{grade}\right]}\ldots\end{multline*}

\begin{multline*}\triangle^\xi v=1\otimes v+v\otimes 1+\sum\xi
(\varepsilon^\mu\otimes\varepsilon^\nu)[e_\mu\otimes(v\wedge e_\nu)-
(v\wedge e_\nu)\otimes e_\mu]+\ldots\end{multline*}

\begin{multline*}\triangle^\xi(v\wedge w)=1\otimes(v\wedge w)+(v\wedge
w)\otimes 1-v\otimes w+w\otimes v\\+\sum\xi(\varepsilon^\mu\otimes
\varepsilon^\nu)[(w\wedge e_\mu)\otimes(e_\nu\wedge v)-
(v\wedge e_\mu)\otimes(e_\nu\wedge w)]+\ldots\end{multline*}
If $\xi=0$ then $1\in\Cl$ is group-like, vectors are $(1,1)$-primitive
and $\triangle^{\xi=0}$ is the Duffin \& Kemmer \& Petiau co-product
[Duffin 1938, Kemmer 1939, 1943].

The above co-product $\triangle^\xi$ (as well as \eqref{c} late on) is
co-unital
\begin{gather}\varepsilon\in\lin(M^\wedge,\Bbbk),\qquad\Bbbk\ni
\varepsilon w=\begin{cases}0&\text{if grade $w\neq 0,$}\\
w&\text{if grade $w=0.$}\end{cases}\label{counit}\end{gather}
However co-unit \eqref{counit} is not an algebra map iff $\eta\neq 0$
and unit $u\in\lin(\Bbbk,M^\wedge)$ is not cogebra map iff $\xi\neq 0,$
\begin{gather*}\eta\neq 0\quad\Longleftrightarrow\quad\varepsilon
\not\in\alg(\wedge^\eta,\Bbbk),\nonumber\\\xi\neq 0
\quad\Longleftrightarrow\quad u\not\in\cog(\Bbbk,\triangle^\xi).
\end{gather*}
\begin{conj} A condition $\xi\circ\eta\neq\id_M$ is a necessary and
sufficient condition that exists an antipode $S\in\End(M^\wedge)$ (see
example below).\end{conj}
If $\xi=0$ then antipode exists and $S|M^{\wedge 2}=0.$

\section{Co-gebra (Co-field) of Co-complex Numbers and Antipode}
In the sequel if $\alpha\in M^*$ then $\alpha^2\in\Bbbk$ stands for
$\xi(\alpha\otimes\alpha)$ and if $v\in M$ then $v^2\in\Bbbk$ stands for
$\eta(v\otimes v).$

In sections 3-4 $\dim_\R M=1.$ Let $i\in M.$ If $\eta(i\otimes i)=-1$
then $\C\simeq\Cl(M,\eta).$

If $\xi\circ\eta=\id_M$ and $\alpha v=1\in\Bbbk,$ then $\alpha^2v^2
\equiv\xi(\alpha\otimes\alpha)\eta(v\otimes v)=1.$

For $j\in\C^*\equiv\lin_\R(\C,\R),$
\begin{gather}\wedge^\eta(1\otimes 1)=1,\quad
\wedge^\eta(i\otimes i)=\wedge^\eta_ii=\eta(i\otimes i)\in\R,\nonumber\\
\wedge^\eta(1\otimes i)=\wedge^\eta_1i=i,\quad
\wedge^\eta(i\otimes 1)=\wedge^\eta_i1=i,\nonumber\\
\triangle^\xi 1=1\otimes 1+\xi(j\otimes j)i\otimes i,\nonumber\\
\triangle^\xi i=1\otimes i+i\otimes 1.\label{c}\end{gather}
For $z\in\C$ and $i^2j^2=1,$ $\triangle z=1\otimes z+
\frac{1}{i^2}i\otimes iz.$
From now on
$$a\equiv i^2j^2=\eta(i\otimes i)\xi(j\otimes j)\in\R.$$
\begin{antipode} With respect to co-unit \eqref{counit}, an antipode
$S\in\lin(\C,\C)$ exists iff $\xi\circ\eta\neq\id_M,$ \ie iff
$a\neq 1,$ and then
\be S\,1=\frac{1}{1-a},\quad S\,i=-\frac{i}{1-a}.\label{anti}
\end{equation}\end{antipode}

\section{Hopf-gebra and Bi-gebra of Complex Numbers} Diagram 4 is a
relation among three tensors, a product $\wedge,$ co-product $\triangle$
and a bin-ary scattering $\sigma.$ The purpose of this section is to
determine from Diagram 4 the set of all possible scatterings
$\sigma\in\End({\C}^{\otimes 2})$ for $(\eta,\xi)$-dependent product and
coproduct \eqref{c} on two-dimensional $\R$-space span by
$\{1,i\}\in\C.$ Then the set $\{\wedge^\eta,\triangle^\xi,\sigma\}$ is a
bin-ary bi-gebra.

If $i^2j^2\neq 1$ then exists the unique scattering
$\sigma\in\End_\R({\C}^{\otimes 2}),$
\begin{align}\sigma(1\otimes 1)&=\left(1-\frac{a^2}{1-a}\right)
1\otimes 1-\frac{j^2}{1-a}i\otimes i,\nonumber\\
\sigma(i\otimes i)&=-\frac{1}{1-a}(i\otimes i+i^2\cdot 1\otimes
1),\nonumber\\\sigma(1\otimes i)&=\frac{1}{1-a}(i\otimes 1+a\cdot
1\otimes i),\nonumber\\\sigma(i\otimes 1)&=\frac{1}{1-a}(1\otimes i+a
\cdot i\otimes 1).\label{scat}\end{align}
The minimum polynomial of \eqref{scat} is of the fourth order
\begin{gather*}b\equiv\frac{1+i^2j^2}{1-i^2j^2},\qquad(\sigma+\id)\circ
(\sigma-b\cdot\id)\circ(\sigma^2+ab\cdot\sigma-b\cdot\id)=0.
\end{gather*}
Therefore $\sigma$ \eqref{scat} is invertible iff $i^2j^2\neq\pm 1.$ We
conjecture that $\sigma\in\End(M^{\otimes 2})$ \eqref{scat} is a
(pre)-braid operator \ie $\sigma$ is a solution of the Artin braid
equation (this indeed is the case if $i^2j^2=0$),
\be(\sigma\otimes{\id}_M)\circ({\id}_M\otimes\sigma)\circ(\sigma\otimes
{\id}_M)=({\id}_M\otimes\sigma)\circ(\sigma\otimes{\id}_M)\circ
({\id}_M\otimes\sigma),\label{artin}\end{equation}
\special{em:linewidth 0.7pt}\linethickness{0.7pt}
$$\begin{picture}(45,23)        
\emline{5}{20}{1}{15}{10}{2}
\emline{15}{10}{3}{15}{5}{4}
\emline{15}{20}{5}{15}{15}{6}
\emline{15}{15}{7}{13}{13}{8}
\emline{12}{12}{9}{8}{8}{10}
\emline{10}{20}{11}{8}{18}{12}
\emline{7}{17.20}{13}{5}{15}{14}
\emline{5}{15}{15}{5}{10}{16}
\emline{5}{10}{17}{10}{5}{18}
\emline{7}{6.88}{19}{5.00}{5.16}{20}
\put(25,13){\makebox(0,0)[cc]{$=$}}
\emline{35}{20}{21}{35}{15}{22}
\emline{35}{15}{23}{45}{5.16}{24}
\emline{40}{20}{25}{45}{15.05}{26}
\emline{45}{15.05}{27}{45}{9.89}{28}
\emline{45}{9.89}{29}{43}{8.17}{30}
\emline{42}{6.88}{31}{40}{5.16}{32}
\emline{45}{20}{33}{43}{18.06}{34}
\emline{42}{17.20}{35}{38}{12.90}{36}
\emline{37}{12}{37}{35}{10}{38}
\emline{35}{9.89}{39}{35}{5.16}{40}\end{picture}$$

\begin{braided} The Clifford bi-gebra
$\{\R^2,\wedge^\eta,\triangle^\xi,\sigma\},$ \eqref{c}-\eqref{scat}, is
$\sigma$-braided iff $\eta=0$ or $\xi=0.$\end{braided}

The formulas \eqref{c}-\eqref{anti}-\eqref{scat} describe two-parameter
$\{i^2,j^2\}$-family of Hopf-gebras (and this include the field of
complex numbers) for which neither unit nor co-unit are respecting
co-product and product respectively.

If $\xi\circ\eta=\id_\C$ \ie if $i^2j^2=1,$ then exists $12$-parameters
family of mappings $\sigma\in\End_\R({\C}^{\otimes 2})$ which fit to
bi-gebra. Among other this include the following solution for
$p+q+r=0\in\R,$
\begin{align*}\sigma(1\otimes 1)&=1\otimes 1,\\
\sigma(1\otimes i)&=i\otimes 1+p\cdot 1\otimes i,\\
\sigma(i\otimes 1)&=1\otimes i+q\cdot i\otimes 1,\\
\sigma(i\otimes i)&=r\cdot i\otimes i-i^2\cdot 1\otimes 1.\end{align*}

\section{Bi-universal Hopf-gebra}
\begin{center}\begin{tabular}{cl}\hline\multicolumn{2}{c}{Some notation}
\\\hline\hline
$\Bbbk$& is a commutative ring\\$\Bbbk$-\textbf{mod}& a
category of $\Bbbk$-modules (of $\Bbbk$-alphabets)\\
$\Bbbk$-\textbf{alg}& a category of associative unital
$\Bbbk$-algebras\\
$T:\,\Bbbk$-\textbf{mod} $\rightarrow\,\Bbbk$-\textbf{alg}& the
tensor algebra functor\\&not additive on morphisms;\\
$F:\,\Bbbk$-\textbf{alg} $\rightarrow\,\Bbbk$-\textbf{mod}& the
forgetful functor;\\
$\otimes$& bifunctor of tensor product:\\
& $\Bbbk$-\textbf{mod}$\,\times\,\Bbbk$-\textbf{mod}$\,\rightarrow
\,\Bbbk$-\textbf{mod}.\\
$\otimes$& means $\otimes_\Bbbk$ if not otherwise stated;\\
$\lin\equiv\lin_\Bbbk,$ $\End\equiv\End_\Bbbk$& are both sided
$\Bbbk$-linear bifunctors;\\$M\in$$\;\Bbbk$-\textbf{mod}&is a
$\Bbbk$-module (a $\Bbbk$-alphabet);\\
$M^\otimes=FTM$&a $\Z$-graded
$\Bbbk$-vocabulary in a $\Bbbk$-alphabet $M$;\\
$M^*\equiv\lin(M,\Bbbk)$&a dual $\Bbbk$-module of co-vectors.\\
\hline\end{tabular}\end{center}

A bi-associative (\ie an associative
and co-associative) and bi-u\-ni\-tal (\ie unital and co-unital)
Hopf-gebra in a braided monoidal category ($\equiv$ a braided Hopf-gebra
or a `braided group') has been introduced by Majid in series of papers
in years 1991-1993. In [Oziewicz et al. 1995] we generalized a braided
Hopf-gebra to pre-braided Hopf-gebra when a braid needs not to be
invertible. This generalization was motivated by the following problem:
does exist pre-braid for which exists a pre-braided
bi-universal (\ie universal and co-universal) Hopf-gebra? This is
illustrated by Diagram 4 in the case of the fixed universal product and
of the co-universal co-product. We showed that pre-braided bi-universal
bi-associative and bi-unital Hopf-gebra exists for zero pre-braid only
[Oziewicz et al. 1995].

For a $\Bbbk$-module $M,$ $M^\otimes$ denotes $\Z$-graded
$\Bbbk$-module (not an algebra) \ie a totality of all finite
sentences in $M.$

By definition functors $T$ and $F$ are adjoint [Kan 1958]: bifunctors
$\lin_\Bbbk(\cdot,F\cdot)$ and $\alg_\Bbbk(T\cdot,\cdot)$ are
naturally equivalent. This means that a natural set bijection holds,
\begin{gather}\mbox{$\forall\;M\times A\in
\Bbbk$-\textbf{mod}$\,\times\,\Bbbk$-\textbf{alg},}\nonumber\\
{\lin}_\Bbbk(M,FA)\;\ni\;\ell\longleftrightarrow\ell^A\;\in\;
\alg_\Bbbk(TM,A),\quad\ell^A|M\equiv\ell.\label{TM}\\
\ell^m\equiv\ell^A=m^\otimes\circ\ell^\otimes\in\alg(TM,A),\nonumber\\
\ell^\triangle\equiv\ell^C=\ell^\otimes\circ\triangle^\otimes\in
\cog(C,ShM),\nonumber\end{gather}

$$\begin{picture}(83,38)      
\put(27,13){\vector(1,-1){0.2}}
\emline{13}{27}{1}{27}{13}{2}
\put(32,23){\vector(0,-1){0.2}}
\emline{32}{29}{3}{32}{23}{4}
\put(32,12){\vector(0,-1){0.2}}
\emline{32}{17}{5}{32}{12}{6}
\put(63,27){\vector(-1,1){0.2}}
\emline{77}{13}{7}{63}{27}{8}
\put(82,17){\vector(0,1){0.2}}
\emline{82}{11}{9}{82}{17}{10}
\put(82,29){\vector(0,1){0.2}}
\emline{82}{23}{11}{82}{29}{12}
\put(62,32){\vector(-1,0){0.2}}
\emline{79}{32}{13}{62}{32}{14}
\put(29,32){\vector(1,0){0.2}}
\emline{12}{32}{15}{29}{32}{16}
\put(8,32){\makebox(0,0)[cc]{$M$}}
\put(32,32){\makebox(0,0)[cc]{$M^\otimes$}}
\put(58,32){\makebox(0,0)[cc]{$M$}}
\put(82,32){\makebox(0,0)[cc]{$M^\otimes$}}
\put(70,33){\makebox(0,0)[cb]{projection}}
\put(20,33){\makebox(0,0)[cb]{injection}}
\put(19,19){\makebox(0,0)[rt]{$\ell$}}
\put(69,19){\makebox(0,0)[rt]{$\ell$}}
\put(32,8){\makebox(0,0)[cc]{$A$}}
\put(32,20){\makebox(0,0)[cc]{$A^\otimes$}}
\put(82,20){\makebox(0,0)[cc]{$C^\otimes$}}
\put(82,8){\makebox(0,0)[cc]{$C$}}
\put(83,14){\makebox(0,0)[lc]{$\triangle^\otimes$}}
\put(83,26){\makebox(0,0)[lc]{$\ell^\otimes$}}
\put(33,26){\makebox(0,0)[lc]{$\ell^\otimes$}}
\put(33,14){\makebox(0,0)[lc]{$m^\otimes$}}
\put(41,4){\makebox(0,0)[ct]{\textsc{Diagram} 7.
The universal and co-universal lifts}}\end{picture}$$

An example of realization of $M$-universal tensor $\Bbbk$-al-gebra is
$TM\simeq\{M^\otimes,\otimes\}$ \ie a $\Z$-graded $\Bbbk$-module
$M^\otimes$ of all finite words in an alphabet $M$ with a concatenation
$\otimes$ as a multiplication of $\grade\,\otimes=0.$

An example of realization of $M$-co-universal co-gebra is
$ShM\simeq\{M^\otimes,\sh\},$ sh is short for `shuffle' [Sweedler 1969],
$T(M^*)\simeq(ShM)^*.$

A deformation of bi-gebra is said to be \textit{preferred} if either
product or co-product is not deformed [Gerstenhaber \& Schack 1992,
Bonneau et al. 1994]. One can consider two pre-braid dependent
bi-associative preferred deformations of $0$-braided bi-universal
Hopf-gebra: an universal Hopf-gebra which is not co-universal (an
universal product is not deformed) and co-universal Hopf-gebra which is
not universal (a co-universal shuffle co-product is not deformed). These
families generalize for arbitrary pre-braid the Sweedler construction
for the switch [Sweedler 1969, chapter XII].

There exists an unique pre-braid dependent homomorphism (extending an
identity mapping on generating space) of universal Hopf-gebra into
co-universal Hopf-gebra [Oziewicz et al. 1995, R\'o\.za\'nski 1996].
This Hopf-gebra homomorphism
\begin{enumerate}
\item is a pre-braid dependent deformation of an identity,
\item commutes with an antipod and
\item for invertible braid coincide with the braid-dependent
`(anti)symmetrizer' introduced by Woronowicz in 1989.\end{enumerate}

The image of this Hopf-gebra homomorphism is said to be an
exterior Hopf-gebra. An exterior Hopf-gebra is co-universal and
pre-braided.

An open problem is to find necessary and sufficient conditions (on
braid, scalar product, Lie algebra, etc.) that exists pre-braided
filtered algebra as tensor-de\-pen\-dent quantizations and (the
Chevalley) deformations of exterior Hopf-gebra.

\end{document}